\documentclass[multi]{cambridge7A}
\usepackage[UKenglish]{babel}
\usepackage[longnamesfirst,sectionbib]{natbib}
\usepackage{chapterbib}
\usepackage{amsmath,amssymb,amsthm,amsfonts}
\usepackage{enumerate,url}

\setcitestyle{numbers,square,comma}% alters natbib behaviour

\allowdisplaybreaks[1]

\providecommand*\Index[1]{#1\index{#1}}
\providecommand*\undex[1]{} % abandoned tag
\begin{document}
\alphafootnotes
\author[R. C. Griffiths and D. Span\'o]
{Robert C. Griffiths\footnotemark\ and Dario Span\'o\footnotemark }
\chapter[Diffusion processes and coalescent trees]%
{Diffusion processes and coalescent trees}
\footnotetext[1]{Department of Statistics, University of Oxford,
  1 South Parks Rd, Oxford OX1 3TG; griff@stats.ox.ac.uk}
\footnotetext[2]{Department of Statistics, University of Warwick,
  Coventry CV4 7AL; D.Spano@warwick.ac.uk}
\arabicfootnotes
\contributor{Robert C. Griffiths
  \affiliation{University of Oxford}}
\contributor{Dario Span\'o
  \affiliation{University of Warwick}}

\renewcommand\thesection{\arabic{section}}
\numberwithin{equation}{section}
\renewcommand\theequation{\thesection.\arabic{equation}}

\begin{abstract}
We dedicate this paper to Sir John Kingman on his 70th Birthday.

In modern mathematical population genetics the ancestral history of a population of genes back in time is described by John Kingman's coalescent tree. Classical and modern approaches model gene frequencies by diffusion processes. This paper, which is partly a review, discusses how coalescent processes are dual to diffusion processes in an analytic and probabilistic sense.

Bochner (1954) and Gasper (1972) were interested in characterizations of processes with Beta stationary distributions and Jacobi polynomial eigenfunctions. We discuss the connection with Wright--Fisher diffusions and the characterization of these processes. Subordinated Wright--Fisher diffusions are of this type. An Inverse Gaussian subordinator is interesting and important in subordinated Wright--Fisher diffusions and is related to the Jacobi Poisson Kernel in orthogonal polynomial theory. A related time-subordinated forest of  non-mutant edges in the Kingman coalescent is novel.
\end{abstract}

\subparagraph{AMS subject classification (MSC2010)}92D25, 60J70, 92D15

\section{Introduction} 
\label{intro}
The Wright--Fisher\index{Wright, S.!Wright--Fisher diffusion} diffusion process $\{X(t),t \geq 0\}$ models the relative frequency of type $a$ genes\index{gene!frequency} in a population with two types of genes $a$ and $A$. Genes are subject to random \Index{drift} and mutation\index{mutation|(} over time. The \Index{generator} of the process is
\begin{equation}
{\cal L} = \frac{1}{2}x(1-x)\frac{\partial^2}{\partial x^2}
+ 
\frac{1}{2}\bigl( -\alpha x + \beta(1-x)\bigr)
\frac{\partial}{\partial x},
\label{generator:0}
\end{equation}
where the mutation rate\index{mutation!rate} $A \to a$ is $\frac{1}{2}\alpha$ and the rate $a \to
A$ is $\frac{1}{2}\beta$. If $\alpha$ and $\beta$ are zero then zero and one
are absorbing states where either $A$ or $a$ becomes fixed in the
population. If $\alpha$, $\beta >0$ then $\{X(t),t\geq 0\}$ is a
reversible\index{reversibility} process with a Beta stationary
density\index{beta distribution|(}
\begin{equation}
f_{\alpha,\beta}(y) = B(\alpha,\beta)^{-1}y^{\alpha - 1}(1-y)^{\beta - 1}, \ 0 < y < 1.
\label{intro:0}
\end{equation}
The transition density\index{transition density|(} has an eigenfunction\index{eigenfunction|(} expansion
\begin{equation}
f(x,y;t) = f_{\alpha,\beta}(y)\biggl \{
1 + \sum_{n=1}^\infty\rho_n^\theta (t)
\widetilde{P}^{(\alpha,\beta)}_n(x)
\widetilde{P}^{(\alpha,\beta)}_n(y)
\biggr \},
\label{intro:1}
\end{equation}
where $\theta = \alpha + \beta$, 
\begin{equation}
\rho_n^\theta(t) = \exp \Bigl \{-\frac{1}{2}n(n+\theta - 1)t\Bigr \},
\label{rhodef}
\end{equation}
 and $\bigl \{\widetilde{P}^{(\alpha,\beta)}_n(y), n \in \mathbb{Z}_+ \bigr \}$
 are orthonormal Jacobi polynomials\index{Jacobi, C. G. J.!Jacobi polynomial|(} on the Beta $(\alpha,\beta)$ distribution, scaled so that 
\begin{equation*}
\mathbb{E}\Bigl [
\widetilde{P}^{(\alpha,\beta)}_m(Y)
\widetilde{P}^{(\alpha,\beta)}_n(Y)
\Bigr ] = \delta_{mn},\ m,n \in \mathbb{Z}_+
\end{equation*}
under the stationary distribution (\ref{intro:0}). The Wright--Fisher diffusion is also known as the Jacobi diffusion\index{Jacobi, C. G. J.!Jacobi diffusion|(} because of the eigenfunction expansion (\ref{intro:1}).
The classical Jacobi polynomials, orthogonal on
\[
(1-x)^\alpha(1+x)^\beta, \ -1 < x < 1,
\]
can be expressed as 
\begin{equation}
P_n^{(\alpha,\beta)}(x) = \frac{(\alpha + 1)_{(n)}}{n!}
\:_2F_1(-n,n+\alpha+\beta+1;\alpha+1;(1-x)/2),
\label{Pdef}
\end{equation}
where $_2F_1$ is a \Index{hypergeometric function}.
The relationship between the two sets of polynomials is that
\[
\widetilde{P}^{(\alpha,\beta)}_n(x) = c_nP_n^{(\beta -1,\alpha -1)}(2x-1),
\]
where
\[
c_n = \sqrt{\frac{(2n+\alpha+\beta-1)(\alpha+\beta)_{(n-1)}n!}
{\alpha_{(n)}\beta_{(n)}}}.
\]
Define
\begin{equation}
\bar{\cal L} = 
\frac{1}{2}
\frac{\partial^2}{\partial x^2}
x(1-x)
- 
\frac{\partial}{\partial x} 
\frac{1}{2}\bigl ( -\alpha x + \beta(1-x)\bigr ),
\label{generator:1}
\end{equation}
the forward generator\index{generator!forward generator} of the process. The Jacobi polynomials are eigenfunctions satisfying, for $n \in \mathbb{Z}_+$,
\begin{eqnarray}
{\cal L}\widetilde{P}^{(\alpha,\beta)}_n(x) &=& -\frac{1}{2}n(n+\theta - 1)\widetilde{P}^{(\alpha,\beta)}_n(x);
\nonumber \\
\bar{\cal L}f_{\alpha,\beta}(x)\widetilde{P}^{(\alpha,\beta)}_n(x) &=& -\frac{1}{2}n(n+\theta - 1)f_{\alpha,\beta}(x)\widetilde{P}^{(\alpha,\beta)}_n(x).
\label{intro:3}
\end{eqnarray}
The well known fact that the Jacobi polynomials\index{Jacobi, C. G. J.!Jacobi polynomial|)} $\bigl \{\widetilde{P}^{(\alpha,\beta)}_n(x)\bigr \}$ satisfy (\ref{intro:3}) implies that they are eigenfunctions\index{eigenfunction|)} with corresponding eigenvalues $\bigl \{\rho^\theta_n(t)\bigr \}$.

In modern mathematical population genetics\index{mathematical genetics} the
ancestral history\index{ancestral history} of a
population back in time is described by John Kingman's elegant coalescent
process \cite{K1982}\index{Kingman, J. F. C.!Kingman coalescent}. The
connection between the coalescent and
Fleming--Viot diffusion
processes\index{Fleming, W. H.!Fleming--Viot diffusion} is made explicit by
Donnelly\index{Donnelly, P. [Donnelly, P. J.]} and Kurtz\index{Kurtz, T. G.}
in \cite{DK1996}, \cite{DK1999} by their
\index{look-down}look-down process. An
approach by Ethier\index{Ethier, S. N.|(} and Griffiths \cite{EG1993} uses
\Index{duality} to show that a 'non-mutant lines of descent'\index{line of descent|(} process
which considers a \Index{forest} of trees\index{tree} back in time to
their first mutations is dual to the Fleming--Viot
infinitely-many-alleles diffusion process. The
\index{allele!two-allele process}two-allele process
$\{X(t), t \geq 0\}$ is recovered from the Fleming--Viot process by a
2-\Index{colouring} of alleles\index{allele} in the infinitely-many-alleles
model. If there is no mutation then the dual process
is the same as the Kingman coalescent
process\index{Kingman, J. F. C.!Kingman coalescent} with an
entrance boundary\index{entrance boundary|(} at infinity. The dual process approach leads to a
transition density expansion in terms of the transition
functions\index{transition function|(} of the
process which counts the number of non-mutant lineages\index{lineage|(} back in
time. It is interesting to make a connection between the eigenfunction
expansion (\ref{intro:1}) and dual process expansion of the transition
densities of $\{X(t),t \geq 0\}$. Bochner \cite{B1954}\index{Bochner, S.|(} and
Gasper \cite{G1972}\index{Gasper, G.} find characterizations of processes
which have Beta stationary distributions\index{beta distribution|)} and Jacobi
polynomial eigenfunctions. Subordinated Jacobi
processes
$\{X\bigl (Z(t)\bigr),t \geq 0 \}$, where $\{Z(t), t \geq 0\}$ is a L{\'e}vy
process\index{Levy, P.@L\'evy, P.!L\'evy process}, fit into this class,
because subordination\index{subordination|(} does not change the eigenvectors or the
stationary distribution of the process.  The subordinated processes are
jump diffusions\index{diffusion!jump diffusion}. A particular class of importance is
when $\{Z(t), t \geq 0\}$ is an Inverse Gaussian
process\index{Gauss, J. C. F.!inverse Gaussian process}. Griffiths
\cite{G2009} obtains characterizations of processes with
stationary distributions in the Meixner\index{Meixner, J.} class, as well as
for Jacobi processes. The current paper is partly a review paper describing
connections between Jacobi
diffusions\index{Jacobi, C. G. J.!Jacobi diffusion|)}, eigenfunction
expansions of transition functions, coalescent trees\index{coalescent tree},
and Bochner\index{Bochner, S.|)} characterizations. Novel results
describe the subordinated\index{subordination|)} non-mutant lines-of-descent
process\index{line of descent|)} when the subordination is with an Inverse
Gaussian process.
\section{A coalescent dual process}
\label{dual}
A second form of the transition density\index{transition density|)} (\ref{intro:1}) derived in Ethier\index{Ethier, S. N.|)} and Griffiths \cite{EG1993} is
\begin{equation}
f(x,y;t) = \sum_{k=0}^\infty q_k^\theta(t)\sum_{l=0}^k{\cal B}(l;k,x)f_{\alpha+l,\beta+k-l}(y),
\label{dual:0}
\end{equation}
where
\[
{\cal B}(l;k,x) = \binom lk x^k(1-x)^{l-k},\ k = 0,1,\ldots ,l
\]
is the Binomial distribution and $\bigl \{q_k^\theta (t)\bigr \}$ are the
transition functions of a \Index{death process} with an entrance boundary of infinity,
and death rates $k(k+\theta -1)/2$, $k \geq 1$. The death process represents
the number of non-mutant ancestral lineages back in time in the coalescent
process with mutation. The number of lineages\index{lineage|)} decreases from $k$ to $k-1$ from
\Index{coalescence} at rate $\binom k2$ or mutation\index{mutation|)} at rate $k\theta/2$.
If there is no mutation, $\{q_k^0(t),t \geq 0\}$ are transition functions of
the number of edges in a Kingman coalescent
tree\index{Kingman, J. F. C.!Kingman coalescent tree}.
There is an explicit expression for the transition functions beginning with
the entrance boundary\index{entrance boundary|)} of infinity
\cite{G1980,T1984,G2006}\index{Tavar\'e, S.} of
\begin{equation}
q_k^{\theta}(t)
=\sum_{j=k}^\infty\rho_j^{\theta}(t)
(-1)^{j-k}
\frac{ (2j+ \theta - 1)(k+\theta)_{(j-1)}}{k!(j-k)!},
\label{trQ}
\end{equation}
recalling that $\rho_n^\theta(t)$ is defined by (\ref{rhodef}).
A complex-variable representation\undex{complex-variable representation|(} of (\ref{trQ}) is found
in \cite{G2006}. Let $\{X_t, t \geq 0\}$ be standard Brownian
motion\index{Brown, R.!Brownian motion}  so $X_t$ is $N(0,t)$. Denote $Z_t=e^{iX_t}$ and $\omega_t = e^{-\frac{1}{2}\theta t}$, then
\begin{equation}
q_k^{\theta}(t) =
e^{\frac{1}{8}t}
\frac{\Gamma (2k+\theta)}{\Gamma (k+\theta )k!}
\mathbb{E}\biggl [
\frac{(\omega_t Z_t)^k(1-\omega_t Z_t)}
{\sqrt{Z_t}(1 + \omega_t Z_t)^{2k+\theta}}\biggr ],
\label{complexrep}
\end{equation}
for $k = 0$, 1, \ldots.
The transition functions for the process beginning at $n$, rather than infinity, are
\begin{equation}
q^\theta_{nk}(t) =\sum_{j=k}^n
\rho^\theta_j(t)
(-1)^{j-k}
\frac{ (2j+ \theta - 1)(k+\theta)_{(j-1)}n_{[j]}}
{k!(j-k)!(n+\theta)_{(j)}},
\label{trQ:1}
\end{equation}
for $k = 0$, 1, \ldots, $n$. An analogous complex-variable representation\undex{complex-variable representation|)} to (\ref{complexrep}) is
\begin{eqnarray}
q^\theta_{nk}(t) &=&
\frac{\Gamma(n+\theta )\Gamma (2k+\theta )}
{\Gamma (k+\theta )\Gamma (n + k + \theta)}
\binom nk e^{\frac{1}{8}(\theta - 1)^2t}
\mathbb{E}\bigl [
Z_t^{k+ (\theta -1)/2}(1-Z_t)
\nonumber \\
&&~~~~~~~~~~~~\times\: _2F_1(-n+k+1,\theta + 2k;n+k+\theta;Z_t)
\bigr ]
\label{complexrep:1}
\end{eqnarray}
for $k=0$, 1, \ldots, $n$.
The expansion (\ref{dual:0}) is derived from a two-dimension\-al
dual\index{duality|(} \Index{death process} $\bigl \{L^\theta
(t) \in \mathbb{Z}^2_+,t \geq 0\bigr \}$ which looks back in time in the
\index{diffusion!diffusion process}diffusion process
$\bigl \{X(t), t \geq 0\bigr \}$. A derivation in
this paper is from \cite{EG2009}\index{Etheridge, A. M.}, which follows more
general analytic derivations in \cite{EG1993}\index{Ethier, S. N.} for a
Fleming--Viot\index{Fleming, W. H.!Fleming--Viot diffusion} model and
\cite{BEG2000}\index{Barbour, A. D.} for a diffusion model with
\Index{selection}.
Etheridge and Griffiths \cite{EG2009} give a very clear probabilistic
derivation in a Moran model\index{Moran, P. A. P.} with selection that provides an understanding of
earlier derivations. A sketch of a derivation of (\ref{dual:0})
from \cite{EG2009} is the following. Let $x_1 = x$, $x_2 = 1-x$ and define for
$k\in \mathbb{Z}_+^2$
\[
g_k(x) = \frac{\theta_{(|k|)}}{{\alpha}_{(k_1)}{\beta}_{(k_2)}}
x_1^{k_1}x_2^{k_2},
\]
then
\begin{equation}
{\cal L} g_k(x) = \frac{1}{2}(|k| + \theta -1)
\bigl [ k_1g_{k-e_1}(x) + k_2g_{k-e_2}(x) - |k|g_k(x)\bigr ].
%\label{introtwod}
\label{dual:1}
\end{equation}
Here and elsewhere we use the notation $|y| = \sum_{j=1}^dy_j$ for a $d$-dimension\-al vector $y$. In this particular case $|k| = k_1 + k_2$.
To obtain a dual process the \Index{generator} is regarded as acting on $k=(k_1,k_2)$, rather than $x$.
The dual process is a two-dimensional death process $\{L^\theta(t), t \geq 0\}$,
the rates of which are read off from the coefficients of the functions 
$g$ on the right-hand side of (\ref{dual:1});
\begin{equation}
k \to k-e_i \text{~~~at rate~~~} \frac{1}{2}\frac{k_i}{|k|}\cdot |k|(|k| + \theta -1).
\label{dual:2}
\end{equation}
The total size, $|L^\theta(t)|$, is a 1-dimensional death process in which
\[
|k| \to |k|-1 \text{~~~at rate~~~} \frac{1}{2}|k|(|k| + \theta -1)
\]
with transition functions denoted by $\bigl \{q_{ml}^{\theta}(t), t \geq 0\bigr \}$.
There is hypergeometric sampling\index{sampling!hypergeometric sampling} of types which do not die, so
\begin{equation}
P\bigl (L(t) = l\bigm| L(0)=m\bigr) = q_{ml}^\theta (t)  
= q^{\theta}_{|m||l|}(t)\frac{\binom{m_1}{l_1}\binom{m_2}{l_2}}{{\binom|m|}{|l|}},
\label{dual:3}
\end{equation}
where $q^\theta_{|m||l|}(t)$ is defined in (\ref{trQ:1}).
The dual\index{duality|)} equation obtained by regarding ${\cal L}$ as acting on $x$ or $k$ in (\ref{dual:1}) is
\begin{equation}
\mathbb{E}_{X(0)}\Bigl [g_{L(0)}\bigl (X(t)\bigr )\Bigr ]
= \mathbb{E}_{L(0)}\Bigl [g_{L(t)}\bigl (X(0)\bigr )\Bigr ],
\label{dual:4}
\end{equation}
where expectation on the left is with respect to the distribution of $X(t)$,
and on the right with respect to the distribution of $L(t)$.  Partitioning the expectation on the right of (\ref{dual:4}) by values taken by $L(t)$,
\begin{eqnarray}
&&\mathbb{E}_{x}\biggl [\binom{|m|}{m_1}X_1(t)^{m_1}X_2(t)^{m_2}\biggr ]
\label{dual:5} \\
&&= \nonumber
\binom{|m|}{m_1}\frac{{\alpha}_{(m_1)}{\beta}_{(m_2)}}{\theta_{(m_1+m_2)}}
\sum_{l \leq m}
x_1^{l_1}x_2^{l_2}\frac{\theta_{(|l|)}}{{\alpha}_{(l_1)}{\beta}_{(l_2)}}
q^{\theta}_{|m||l|}(t)\binom{|l|}{l_1}\frac{{m_1}_{[l_1]}{m_2}_{[l_2]}}{|m|_{[|l|]}}.
\end{eqnarray}
The transition distribution of $X(t)$ now has an expansion derived from 
an inversion formula applied to (\ref{dual:5}).
Letting
$m_1$, $m_2 \to \infty$ with $m_1/|m| \to y_1$, $m_2/|m| \to y_2$ gives
\[
f(x,y;t) = \sum_{l \in \mathbb{Z}_+^2}
q^{\theta}_{|l|}(t)
\binom{|l|}{l_1}x_1^{l_1}x_2^{l_2}B(\alpha+l_1,\beta+l_2)^{-1}y_1^{l_1+\alpha-1}y_2^{l_2+\beta-1},
\]
which is identical to (\ref{dual:0}).

The two-allele\index{allele!two-allele process} Wright--Fisher
diffusion\index{Wright, S.!Wright--Fisher diffusion|(} is a special case of a
much more general Fleming--Viot\index{Fleming, W. H.!Fleming--Viot diffusion}
measure-valued\index{diffusion!measure-valued diffusion} diffusion process which has
${\cal P}(S)$, the probability measures on $S$, a compact metric space, as a
state space.  The mutation operator\index{mutation!operator} in the process is
\[
(Af)(x) = \frac{\theta}{2}\int_S \bigl (f(\xi ) - f(x)\bigr )\nu_0(d\xi ),
\]
where $\nu_0 \in {\cal P}(S)$ and $f:S \to \mathbb{R}$. The \Index{stationary
measure} is a
Poisson--Dirichlet\index{Poisson, S. D.!PoissonDirichlet random measure@Poisson--Dirichlet random measure|(}
(Ferguson--Dirichlet)\index{Ferguson, T. S.!Ferguson--Dirichlet random measure}
random measure
\[
\mu = \sum_{i=1}^\infty x_i \delta_{\xi_i},
\]
where $\{x_i\}$ is a Poisson--Dirichlet point
process,\index{Poisson, S. D.!PoissonDirichlet point process@Poisson--Dirichlet point process}
${\cal PD}(\theta$), independent of $\{\xi_j\}$ which are 
{\it i.i.d.} $\nu_0 \in {\cal P}(S)$.
A description of the ${\cal PD}(\theta)$ distribution is contained in
Kingman \cite{K1993}.\index{Kingman, J. F. C.}

Denote the stationary distribution of the \Index{random measure} as
\[
\Pi_{\theta,\nu_0}(\cdot ) = \mathbb{P}(\mu \in \cdot).
\]
Ethier\index{Ethier, S. N.} and Griffiths \cite{EG1993} derive a transition function\index{transition function|)} expansion for
$P(t,\mu,\allowbreak d\nu)$ with given initial $\mu \in {\cal P}(S)$ of
\begin{eqnarray}
\mathbb{P}(t,\mu,.) &=& q_0^\theta(t)\Pi_{\theta,\nu_0}(\cdot )\cr
&&~~+ \sum_{n=1}^\infty q_n^\theta(t)
\int_{S^n}\mu^n(dx_1\times \cdots \times dx_n)
\cr
&&~~~~~~~~~~
\times\:\Pi_{n+\theta,(n+\theta)^{-1}\{n\eta_n(x_1,\ldots ,x_n)+ \theta\nu_0\}}(\cdot ),
\label{dual:6}
\end{eqnarray}
where $\eta_n(x_1,\ldots ,x_n)$ is the \Index{empirical measure} of points
$x_1$, \ldots, $x_n \in S$:
\[
\eta_n(x_1,\ldots ,x_n) = n^{-1}(\delta_{x_1} + \cdots + \delta_{x_n}).
\]
There is the famous Kingman coalescent process
tree \cite{K1982}\index{Kingman, J. F. C.!Kingman coalescent tree} behind the
pretty representation (\ref{dual:6}). The coalescent tree has an
\Index{entrance boundary} at infinity and a coalescence rate of $\binom k2$
while there are $k$ edges in the tree. Mutations\index{mutation|(} occur
according to a Poisson process\index{Poisson, S. D.!Poisson process} of rate
$\theta/2$ along the edges of the coalescent
tree. $\bigl \{q_n^\theta(t)\bigr \}$ is the distribution of the number of
non-mutant edges in the tree at time $t$ back. The number of non-mutant edges
is the same as the number of edges in a \Index{forest}
where \Index{coalescence} occurs to non-mutant edges and trees are rooted back
in time when mutations occur on an edge. If the time origin is at time $t$
back and there are $n$ non-mutant edges at the origin then the leaves of the
infinite-leaf tree represent the population at $t$ forward in time divided
into relative frequencies of families of types which are either the $n$
non-mutant types chosen at random from time zero, or mutant types chosen from
$\nu_0$ in $(0,t)$. The frequencies of non-mutant families, scaled to have a
total frequency unity, have a Dirichlet
distribution\index{Dirichlet, J. P. G. L.!Dirichlet distribution|(} with unit
index parameters, and the new mutation families, scaled to have total
frequency unity, are distributed according to a Poisson--Dirichlet random
measure\index{Poisson, S. D.!PoissonDirichlet random measure@Poisson--Dirichlet random measure|)} with rate
$\theta$ and type measure $\nu_0$. The total frequency of new mutations has a
Beta $(\theta,n-1)$ distribution\index{beta distribution|(}. An extended
description of the tree process\index{tree!process} is in Griffiths \cite{G2006}.

A $d$-dimensional reversible\index{reversibility}
\index{diffusion!diffusion process}diffusion process
model for gene frequencies\index{gene!frequency} which arises as a limit from
the Wright--Fisher model\index{Wright, S.!Wright--Fisher diffusion|)} has a
backward generator\index{generator!backward generator}
\begin{equation}
{\cal L}=
\frac{1}{2}
\sum_{i=1}^d\sum_{j=1}^d
x_i(\delta_{ij} - x_j)\frac{\partial^2}{\partial x_i\partial x_j} 
+
\frac{1}{2}
\sum_{i=1}^d(\epsilon_i - \theta x_i)\frac{\partial}{\partial x_i},
\label{Generator_b}
\end{equation}
where $\theta = |\epsilon|$.
In this model mutation is parent-independent from type $i \to j$ at rate
$\frac{1}{2}\epsilon_j$, $i$, $j = 1$, \ldots, $d$.
Assuming that $\epsilon > 0$, the stationary density is the Dirichlet
density\index{Dirichlet, J. P. G. L.!Dirichlet distribution|)}
\begin{equation}
%{\cal D}(x,\epsilon) = 
\frac{\Gamma (\theta)}{\Gamma (\epsilon_1) \cdots \Gamma (\epsilon_d)}
x_1^{\epsilon_1-1}\cdots x_d^{\epsilon_d-1},
\label{Dirichlet}
\end{equation}
for $x_1$, \ldots, $x_d > 0$ and $\sum_1^d x_i = 1$.
Griffiths \cite{G1979} shows that the transition
density\index{transition density|(} in the model has
eigenvalues 
\[
\rho_{|n|}(t) = e^{-\frac{1}{2}|n|(|n|+\theta-1)t}
\]
repeated 
\[
\binom{|n| + d - 2}{|n|}
\]
times corresponding to eigenvectors $\bigl \{Q_n^\circ (x),
n \in \mathbb{Z}_+^{d-1}\bigr \}$ which are multitype orthonormal\index{orthonormality|(} polynomials
of total degree $|n|$ in $x$. As eigenfunctions\index{eigenfunction} the polynomials satisfy
\begin{equation}
{\cal L}Q_n^\circ (x) = - \frac{1}{2}|n|(|n|+\theta-1)Q_n^\circ (x).
\label{backeq}
\end{equation}
The eigenvalues $\{\rho_k(t), k \in \mathbb{Z}_+\}$ do not depend on the dimension $d$.
The transition density with $X(0) = x$, $X(t) = y$ has the form
\begin{equation}
f(x,y,t)=
{\cal D}(y,\epsilon)\biggl \{1 + \sum_{|n|=1}^\infty \rho_{|n|}(t)Q_{|n|}(x,y)\biggr \}.
\label{densityone}
\end{equation}
The kernel polynomials\index{kernel polynomial} on the Dirichlet $\{Q_{|n|}(x,y)\}$ appearing in (\ref{densityone})
are defined as
\begin{equation}
Q_{|n|}(x,y) 
= \sum_{\{n:|n|{\rm ~fixed}\}}Q^\circ_n(x)Q^\circ_n(y)
\end{equation}
for \emph{any} complete orthonormal polynomial set $\{Q^\circ_n(x)\}$ on the Dirichlet distribution (\ref{Dirichlet}). 
If $d=2$,
\[
Q_{|n|}(x,y) =
 \widetilde{P}^{(\epsilon_1,\epsilon_2)}_{|n|}(x)
 \widetilde{P}^{(\epsilon_1,\epsilon_2)}_{|n|}(y)
\]
where $\bigl \{ \widetilde{P}^{(\epsilon_1,\epsilon_2)}_{|n|}(x)\bigr \}$ are
orthonormal Jacobi polynomials\index{Jacobi, C. G. J.!Jacobi polynomial} on
the Beta distribution\index{beta distribution|)} on $[0,1]$.  In
general $n$ is just a convenient index system for the polynomials since the
number of polynomials of total degree $|n|$ is always the same as the number
of solutions of $n_1+\cdots + n_{d-1} = |n|$,
\[
\binom{|n| + d -2}{|n|}.
\]
$Q_{|n|}(x,y)$ is invariant under the choice of which orthonormal\index{orthonormality|)} polynomial set is used.
The individual polynomials $Q_n^\circ(x)$ are uniquely determined by their leading coefficients of degree $|n|$ and $Q_{|n|}(x,y)$. A specific form is
\begin{equation}
Q_{|n|}(x,y) = (\theta + 2|n| -1)
\sum_{m=0}^{|n|}(-1)^{|n|-m}
\frac{(\theta + m)_{(|n|-1)}}{m!(|n|-m)!}\xi_m,
\label{Dpoly}
\end{equation}
where 
\begin{equation}
\xi_m =
\sum_{|l|=m}\binom ml
 \frac{ \theta_{(m)} }
{ \prod_1^d{\epsilon_i}_{(l_i)} }
\prod_1^d(x_iy_i)^{l_i}.
\end{equation}
An inverse relationship is
\begin{equation}
\xi_m = 1 + \sum_{|n|=1}^m 
\frac{m_{[|n|]}}{(\theta + m)_{(|n|)}}
Q_{|n|}(x,y).
\label{Inverse}
\end{equation}
The transition distribution (\ref{densityone}) is still valid if any or all elements of $\epsilon$ are zero. 
The constant term in the expansion then vanishes as the
\index{diffusion!diffusion process}diffusion process is
transient\index{transience} and there is not a stationary distribution.
For example, if $\epsilon=0$, 
\begin{equation}
f(x,y,t)=
\prod_{j=1}^dy_j^{-1}\Biggl \{\sum_{|n| \geq d}^\infty \rho_{|n|}(t)Q^0_{|n|}(x,y)\Biggr \},
\label{densityzero}
\end{equation}
where
\begin{equation}
Q^0_{|n|}(x,y) = (2|n| -1)
\sum_{m=1}^n(-1)^{|n|-m}
\frac{(m)_{(|n|-1)}}{m!(|n|-m)!}\xi^0_m,
\end{equation}
with
\begin{equation}
\xi^0_m =
\sum_{\{l: l > 0,|l|=m\}}\binom ml
 \frac{(m-1)!}
{\prod_1^d(l_i-1)!}
\prod_1^d(x_iy_i)^{l_i}.
\end{equation}
The derivation of (\ref{densityone}) is a very classical approach.  The same
process can be thought of as arising from an infinite-leaf
\Index{coalescent tree} similar to the description in the
Fleming--Viot infinitely-many-alleles
process. The coalescent rate while there are $k$ edges in the tree is $\binom
k2$ and mutations occur along edges at rate $\theta/2$. In this model there
are $d$ types, 1, 2, \ldots, $d$ and the probability of mutation $i \to j$,
given a mutation, is $\epsilon_j/\theta$. This is equivalent to a
$d$-colouring of alleles\index{allele} in the
Fleming--Viot\index{Fleming, W. H.!Fleming--Viot diffusion}
infinitely-many-alleles
model. Think backwards from time $t$ back to time 0. Let $y = (y_1,\ldots
,y_d)$ be the relative frequencies of types in the infinite number of leaves
at the current time $t$ forward and $x = (x_1,\ldots ,x_d)$ be the frequencies
in the population at time 0. Let $l$ be the number of non-mutant edges at time
$0$ which have families at time $t$ in the leaves of the tree. Given these $l$
edges let $U = (U_1,\ldots ,U_l)$ be their relative family sizes in the
leaves, and $V = (V_1,\ldots ,V_d)$ be the frequencies of families derived
from new mutations on the tree edges in $(0,t)$. The distribution of $U \oplus
V = (U_1,\ldots ,U_l,V_1,\ldots ,V_d)$ is ${\cal D}(u\oplus v, (1,\ldots
,1) \oplus \epsilon)$. The type of the $l$ lines, and therefore their
families, is chosen at random from the frequencies $x$.  The distribution of
the number of non-mutant lines at time $0$ from the population at $t$ is
$q^\theta_l(t)$.  The transition density in the diffusion (\ref{densityone})
is identical to the mixture distribution arising from the coalescent\index{coalescence}
\begin{equation}
f(x,y,t) =
\sum_{|l|=0}^\infty q_{|l|}^{\theta}(t)
\sum_{\{l: |l|{\rm ~fixed}\}}{\cal M}(l,x){\cal D}(y,\epsilon + l),
\label{densitytwo}
\end{equation}
by considering types of non-mutant lines, and adding Dirichlet
variables\index{Dirichlet, J. P. G. L.!Dirichlet distribution} and
parameters according to $l_i$ non-mutant families being of type $i$.  ${\cal
M}(l,x)$ is the multinomial distribution describing the choice of the initial
line types from the population at time 0.  The expansion when $d=2$
corresponds to (\ref{intro:1}).  The argument is valid if any elements of
$\epsilon$ are zero, considering a generalized Dirichlet
distribution\index{Dirichlet, J. P. G. L.!generalized Dirichlet distribution}
${\cal D}(x,\epsilon)$ where if $\epsilon_i = 0$, then $X_i = 0$ with
probability 1.

The algebraic identity of (\ref{densitytwo}) and (\ref{densityone}) is easy to see by expressing $Q_{|n|}(x,y)$ in terms of $\{\xi_m\}$, then collecting coefficients of $\xi_{|l|}$ in (\ref{densityone}) to obtain (\ref{densitytwo}).
Setting $\rho_0(t) = 1$ and $Q_0(x,y) =1$, the transition density is
\begin{eqnarray}
f(x,y,t)  &=&
{\cal D}(y,\epsilon)\sum_{|n|=0}^\infty\rho_{|n|}(t)Q_{|n|}(x,y)
\nonumber \\
&=&\sum_{l\in \mathbb{Z}_+^d}^\infty\left [
\sum_{|n|=|l|}^\infty\rho_{|n|}(t)(\theta+2|n|-1)(-1)^{|n|-|l|}
\frac{(\theta + |l|)_{(|n|-1)}}{|l|!(|n|-|l|)!}\right ]
\nonumber \\
&&~~~~~~~~~~~~~~~~~~~~~~~~~~~~~~~~~~~~~~~~~~~~~
\times\:{\cal D}(y,\epsilon)\xi_l(x,y)
\nonumber \\
&=&\sum_{|l|=0}^\infty q_{|l|}^\theta (t)
\sum_{\{l: |l|{\rm ~fixed}\}}{\cal M}(l,x){\cal D}(y,\epsilon + l).
\label{OFLD}
\end{eqnarray}
The non-mutant line-of-descent\index{line of descent} process with transition
probabilities\break $\{q_n^\theta(t)\}$ appears in all the Wright--Fisher
diffusion\index{Wright, S.!Wright--Fisher diffusion|(} processes mentioned in
this section as a fundamental dual\index{duality} process. The process does
not depend on the dimension of the diffusion, partly because the
$d$-dimensional process can be recovered from the
\index{measure valued process@measure-valued process}measure-valued process
as a special case by \Index{colouring} new
mutations\index{mutation|)} into
$d$ classes with probabilities $(\epsilon_1/\theta,\epsilon_2/\theta,\ldots
,\epsilon_d/\theta)$ with $\theta=\sum_{j=1}^d\epsilon_j$. It is also
interesting to see the derivation of the $d$-dimensional transition
density\index{transition density|)} expansion as a mixture in terms of
$\{q_n^\theta(t)\}$ via the \Index{orthogonal-function expansion} of the
transition density in (\ref{OFLD}).
\section{Processes with beta stationary distributions and Jacobi polynomial eigenfunctions} 
In this section we consider 1-dimensional processes which have Beta stationary
distributions\index{beta distribution!bivariate beta distribution|(} and Jacobi
polynomial\index{Jacobi, C. G. J.!Jacobi polynomial}
eigenfunctions\index{eigenfunction}, and their connection with
Wright--Fisher diffusion\index{Wright, S.!Wright--Fisher diffusion|)}
processes. We begin by considering Bochner \cite{B1954}\index{Bochner, S.} and
Gasper's\index{Gasper, G.|(} \cite{G1972} characterization of bivariate Beta distributions.

A class of bivariate distributions with Beta marginals and Jacobi polynomial eigenfunctions has the form 
\begin{equation}
f(x,y) = 
f_{\alpha\beta}(x)
f_{\alpha\beta}(y)
\biggl \{
1 + \sum_{n=1}^\infty\rho_n
\widetilde{P}^{(\alpha,\beta)}_n(x)
\widetilde{P}^{(\alpha,\beta)}_n(y)
\biggr \},
\label{general:1}
\end{equation}
where $\{\rho_n, n \in \mathbb{Z}_+\}$ is called a
correlation sequence\index{correlation!sequence|(}. The transition density
(\ref{intro:1}) in the
Jacobi diffusion\index{Jacobi, C. G. J.!Jacobi diffusion} has
the form of the conditional density of $Y$ given $X=x$ in (\ref{general:1})
with $\rho_n \equiv \rho_n^\theta (t)$. Bochner \cite{B1954} and
Gasper \cite{G1972} worked on characterizations of sequences $\{\rho_n\}$ such
that the expansion (\ref{general:1}) is positive, and thus a probability
distribution. It is convenient to normalize the Jacobi polynomials by taking
\[
R_n^{(\alpha,\beta)}(x) = 
\frac{\widetilde{P}_n^{(\alpha,\beta)}(x)}{\widetilde{P}_n^{(\alpha,\beta)}(1)}
\]
so that $R_n^{(\alpha,\beta)}(1) =  1$;  denote
\[
h^{-1}_n = \mathbb{E}\bigl [R_n^{(\alpha,\beta)}(X)^2\bigr ]
= \frac{(2n + \alpha+\beta - 1)(\alpha + \beta)_{(n-1)}\beta_{(n)}}
{\alpha_{(n)}n!},
\]
and write
\begin{equation}
f(x,y) = 
f_{\alpha\beta}(x)
f_{\alpha\beta}(y)
\biggl \{
1 + \sum_{n=1}^\infty\rho_nh_n
R^{(\alpha,\beta)}_n(x)
R^{(\alpha,\beta)}_n(y)
\biggr \}.
\label{general:2}
\end{equation}
Bochner \cite{B1954} defined a bounded sequence $\{c_n\}$ to be positive
definite\index{positive definite|(}
with respect to the Jacobi polynomials if
\[
\sum_{n\geq 0} a_nh_nR^{(\alpha ,\beta)}_n(x) \geq 0,\>\> \sum_{n\geq 0}|a_n|h_n < \infty
\]
implies that
\[
\sum_{n\geq 0} a_nc_nh_nR^{(\alpha ,\beta)}_n(x) \geq 0.
\] 
Then $\{\rho_n\}$ is a correlation
sequence \emph{if and only if} it is a positive
definite\index{positive definite|)} sequence.  The \emph{only if} proof
follows from
\[
\sum_{n\geq 0}a_n\rho_nR_n^{(\alpha,\beta)}(x)
= \mathbb{E}\bigg [\sum_{n\geq 0}a_nh_n
R_n^{(\alpha,\beta)}(Y) \biggm| X = x\bigg ] \geq 0,
\]
where $(X,Y)$ has the distribution (\ref{general:2}).
The \emph{if} proof follows at least heuristically by noting that
\[
\sum_{n\geq 0}h_n R_n^{(\alpha,\beta)}(x)R_n^{(\alpha,\beta)}(y)
= \frac{\delta (x-y)}{f_{\alpha,\beta}(x)} \geq 0,
\]
where $\delta (\cdot)$ has a unit point mass at zero,
so if $\{\rho_n\}$ is a positive definite sequence then
\[
\sum_{n\geq 0}\rho_nh_n R_n^{(\alpha,\beta)}(x)R_n^{(\alpha,\beta)}(y) \geq 0
\]
and (\ref{general:2}) is non-negative. A careful proof is given in \cite{G1970}.

Under the conditions  that 
\begin{equation}
\alpha < \beta \text{~and~either~}1/2 \leq \alpha\text{~or~}\alpha + \beta \geq 2,
\label{XX}
\end{equation}
it is shown in \cite{G1972}\index{Gasper, G.|)} that a sequence $\rho_n$ is positive definite if and only if
\begin{equation}
\rho_n = \mathbb{E}\bigl [R_n^{(\alpha ,\beta )} (Z)\bigr ]
\label{XX:0}
\end{equation}
for some random variable $Z$ in $[0,1]$.
If the conditions (\ref{XX}) do not hold then there exist $x$, $y$, $z$ such that $K(x,y,z) < 0$.
The sufficiency rests on showing that under the conditions (\ref{XX}) for $x$,
$y$, $z \in [0,1]$,
\begin{equation}
K(x,y,z) = \sum_{n=0}^\infty h_n
R_n^{(\alpha ,\beta )} (x)
R_n^{(\alpha ,\beta )} (y)
R_n^{(\alpha ,\beta )} (z) \geq 0.
\label{K:3}
\end{equation}
The sufficiency of (\ref{XX:0}) is then clear by mixing over a distribution for $Z$ in (\ref{K:3}) to get positivity.
The necessity follows by setting $x=1$ in 
\[
\rho_n R_n^{(\alpha ,\beta )} (x) = \mathbb{E}\bigl [ R_n^{(\alpha ,\beta )} (Y) \bigm| X=x\bigr ],
\]
and recalling that $R_n^{(\alpha,\beta)}(1) = 1$, so that $Z$ is distributed
as $Y$ conditional on $X=1$.  This implies that extreme correlation
sequences\index{correlation!sequence} in exchangeable bivariate Beta
distributions\index{beta distribution!bivariate beta distribution|)} with Jacobi polynomial
eigenfunctions\index{eigenfunction} are the scaled Jacobi
polynomials\index{Jacobi, C. G. J.!Jacobi polynomial|(} $\bigl \{R^{(\alpha ,\beta
)}_n(z), z \in [0,1] \bigr \}$.  Bochner \cite{B1954}\index{Bochner, S.} was
the original author
to consider such problems for the ultraspherical
polynomials\index{ultraspherical polynomial}, essentially
orthogonal polynomials on Beta distributions with equal parameters.

A characterization of reversible\index{reversibility} Markov
processes\index{Markov, A. A.!Markov process|(} with stationary Beta
distribution and Jacobi polynomial eigenfunctions\index{eigenfunction}, from \cite{G1972}, under\break (\ref{XX}), is that they have transition functions of the form 
\begin{equation}
f(x,y;t) = 
f_{\alpha\beta}(y)
\biggl \{
1 + \sum_{n=1}^\infty c_n(t)h_n
R^{(\alpha,\beta)}_n(x)
R^{(\alpha,\beta)}_n(y)
\biggr \},
\label{general:4}
\end{equation}
with $c_n(t) = \exp\{-d_nt\}$, where
\begin{equation}
d_n = \sigma n(n + \alpha + \beta  -1) + 
\int_0^{1-}\>\frac{1 - R_n^{(\alpha ,\beta )}(z)}{1-z}\>\nu(dz),
\label{ezero}
\end{equation}
$\sigma \geq 0$, and $\nu$ is a finite measure on $[0,1)$.
If $\nu (\cdot) \equiv 0$, a null measure, then $f(x,y;t)$ is the transition
function of a Jacobi diffusion\index{Jacobi, C. G. J.!Jacobi diffusion|(}.
%Characterizations of correlation sequences (\ref{XX:0}) and (\ref{ezero}) are not complete for all parameters $\alpha, \beta >0$. They are both necessary conditions for all $\alpha,\beta >0$, but only sufficient if (\ref{XX}) holds.

Eigenvalues of a general reversible time-homogeneous Markov process with
countable spectrum must satisfy Bochner's\index{Bochner, S.} consistency conditions:
\begin{enumerate}
\item[(i)] $\{c_n(t)\}$ is a correlation sequence for each $t\geq 0$,
\item[(ii)] $c_n(t)$ is continuous in $t \geq 0$,
\item[(iii)] $c_n(0)=c_0(t)=1$, and
\item[(iv)] $c_n(t+s) = c_n(t)c_n(s)$ for $t$, $s \geq 0$.
\end{enumerate}
If there is a \Index{spectrum} $\{c_n(t)\}$ with corresponding eigenfunctions $\{\xi_n\}$ then
\begin{eqnarray*}
c_n(t+s)\xi_n\bigl (X(0)\bigr ) &=& 
\mathbb{E}\Bigl [ \xi_n\bigl (X(t+s)\bigr ) \Bigm| X(0) \Bigr ] \\
&=&
\mathbb{E}\Bigl [ \mathbb{E}\bigl [\xi_n\bigl (X(t+s)\bigr )\bigm| X(s)\bigr ]\Bigm| X(0)\Bigr ] \\
&=&
c_n(t)
\mathbb{E}\Bigl [ \xi_n\bigl (X(s)\bigl ) \Bigm| X(0)\Bigr ] \\
&=&
c_n(t)c_n(s) \xi_n \bigl (X(0)\bigr ),
\end{eqnarray*}
showing (iv).
If a stationary distribution exists and $X(0)$ has this distribution then
the eigenfunctions can be scaled to be orthonormal\index{orthonormality} on
this distribution and the eigenfunction property is then
\[
\mathbb{E}\Bigl [ \xi_m\bigl (X(t)\bigr ) \xi_n\bigl (X(0)\bigr )\Bigr ] = c_n(t)\delta_{mn}.
\]
$\{X(t),t\geq 0\}$ is a Markov process\index{Markov, A. A.!Markov process|)}
such that the transition distribution\index{transition density} of $Y = X(t)$
given $X(0)=x$ is
\begin{equation}
f(x,y;t) = f(y)
\biggl \{1 + \sum_{n=1}^\infty c_n(t)\xi_n(x)\xi_n(y)\biggr \},
\label{Gammatransition}
\end{equation}
where $f(y)$ is the stationary distribution.  In our context $\{\xi_n\}$ are
the orthonormal Jacobi
polynomials.  A Jacobi
process
$\{X(t),t\geq 0\}$ with transition
distributions (\ref{general:4}) can be constructed in the following way, which
is analogous to constructing a general L{\'e}vy
process\index{Levy, P.@L\'evy, P.!L\'evy process} from a compound
Poisson process\index{Poisson, S. D.!compound Poisson process}. Let
$\{X_k,k\in \mathbb{Z}_+ \}$ be a Markov
chain\index{Markov, A. A.!Markov chain} with
stationary distribution $f_{\alpha\beta}(y)$ and transition distribution of
$Y$ given $X=x$ corresponding to (\ref{general:1}), with (\ref{XX}) holding,
and $\{N(t),t \geq 0\}$ be an independent Poisson
process\index{Poisson, S. D.!Poisson process} of rate $\lambda$.
Then $(X_0,X_k)$ has a
\index{correlation!sequence}correlation sequence $\{\rho_n^k \}$ and the
transition functions of $X(t) = X_{N(t)}$ have the form
(\ref{Gammatransition}), with
\begin{equation}
d_n = \lambda
\int_0^1\bigl (1-R_n^{(\alpha,\beta)}(z)\bigr )\mu(dz),
\label{Poissondn}
\end{equation}
where $\mu$ is a probability measure on $[0,1]$.
The general form (\ref{ezero}) is obtained by choosing a pair 
$(\lambda, \mu_\lambda)$ such that
\begin{equation}
d_n = \lim_{\lambda \to \infty}\lambda
\int_0^1\bigl (1-R_n^{(\alpha,\beta)}(z)\bigr )\mu_\lambda(dz)
= \int_0^1 \frac{1-R^{(\alpha,\beta)}_n(z)}{1-z}\nu(dz).
\label{pairlimit}
\end{equation}
Equation (\ref{pairlimit}) agrees with (\ref{ezero}) when any atom $\nu(\{1\})$ is taken out of the integral, because
\[
\lim_{z\to 1}
 \frac{1-R^{(\alpha,\beta)}_n(z)}{1-z} = cn(n+\theta - 1),
\]
where $c \geq 0$ is a constant.
\section{Subordinated Jacobi diffusion processes}
Let $\{X(t), t \geq 0\}$ be a process with transition functions
(\ref{general:4}), and $\{Z(t), t \geq 0\}$ be a non-negative L{\'e}vy process
with Laplace transform\index{Laplace, P.-S.!Laplace transform}
\begin{equation}
\mathbb{E}\bigl [e^{-\lambda Z(t)}\bigr ]
= \exp\left \{-t\int_0^\infty \frac{1 - e^{-\lambda y}}{y}H(dy)\right \},
%\label{Levyrep}
\end{equation}
where $\lambda \geq 0$ and $H$ is a finite measure. 
The subordinated\index{subordination|(} process $\{\widetilde{X}(t)=X(Z(t)),
t \geq 0\}$ is a Markov process which belongs to the same class of processes
with correlation sequences
\begin{equation}
\widetilde{c}_n(t) = \mathbb{E}\bigl [c_n\bigl (Z(t)\bigr )\bigr ]
= \exp\left \{-t\int_0^\infty \frac{1 - e^{-d_ny}}{y}H(dy)\right \},
\label{Levyrep}
\end{equation}
where $H$ is a finite measure. 
$\widetilde{c}_n(t)$ necessarily has a representation as $e^{-\widetilde{d}_nt}$, where $\widetilde{d}_n$ has the form (\ref{pairlimit}) for some measure $\widetilde{\nu}$. We describe the easiest case from which the general case can be obtained as a limit. Suppose
\[
\lambda = \int_0^\infty\frac{H(dy)}{y} < \infty,
\]
and write
\[
G(dy) = \frac{H(dy)}{\lambda y},
\]
so that $G$ is a probability measure. Let
\[
K(dz) = f_{\alpha\beta}(z)\,dz\biggl \{1 + \sum_{n=1}^\infty
h_nR^{(\alpha,\beta)}_n(z)\int_0^\infty e^{-d_ny}G(dy)\biggr \}.
\]
Then $K$ is a probability measure and
\[
\lambda \int_0^1\Bigl (1 - R_n^{(\alpha,\beta)}(z) \Bigr )K(dz)
= \lambda \int_0^\infty \Bigl (1 - e^{-d_ny}\Bigr ) G(dy).
\]
The representation (\ref{pairlimit}) is now obtained by setting
\[
\widetilde{\nu}(dz) = \lambda (1-z)K(dz).
\]
We now consider subordinated Jacobi diffusion
processes\index{Jacobi, C. G. J.!Jacobi diffusion|)}. The
subordinated\index{subordination|)} process is no longer a diffusion process
because $\{Z(t),t\geq 0\}$ is a jump process\index{jump process} and therefore
$\{\widetilde{X}(t),t \geq 0\}$ has discontinuous sample paths.  It is
possible to construct processes such that (\ref{Levyrep}) holds with $d_n = n$
by showing that $e^{-tn}$ is a correlation
sequence\index{correlation!sequence|)} and thus so is
$\mathbb{E}\bigl [e^{-Z(t)n}\bigr ]$.  The construction follows an idea
in \cite{B1954}\index{Bochner, S.|(}.  The Jacobi--Poisson
kernel\index{Jacobi, C. G. J.!JacobiPoisson kernel@Jacobi--Poisson kernel|(} in orthogonal
polynomial theory is
\begin{equation}
1 + \sum_{n=1}^\infty r^nh_n
R^{(\alpha ,\beta)}_n(x)R_n^{(\alpha ,\beta)}(y),
\label{PoissonK}
\end{equation}
which is non-negative for all $\alpha$, $\beta >0$, $x$, $y\in [0,1]$, and
$0 \leq r \leq 1$, for which see
\cite{AAR1999}\index{Andrews, G. E.}\index{Askey, R,}\index{RoyR@Roy, R.}, p112.
The series (\ref{PoissonK}) is a classical one evaluated early in research on
Jacobi polynomials (see \cite{B1938}). In terms of the original Jacobi
polynomials\index{Jacobi, C. G. J.!Jacobi polynomial|)},
(\ref{Pdef}) 
\begin{equation}
\begin{split}
&\sum_{n=0}^\infty r^n\phi_nP_n^{(\alpha,\beta)}(x)P_n^{(\alpha,\beta)}(y)\\
&= \frac{\Gamma(\alpha+\beta+2)(1-r)}
{2^{\alpha+\beta+1}\Gamma (\alpha+1)\Gamma (\beta + 1)(1+r)^{\alpha+\beta+2}}\\
&\times \sum_{m,n=0}^\infty
\frac{
\bigl ((\alpha+\beta+2)/2\bigr )_{(m+n)}
\bigl ((\alpha+\beta+3)/2\bigr )_{(m+n)}
}
{
(\alpha+1)_{(m)}(\beta+1)_{(m)}m!n!
}
\biggl(\frac{a^2}{k^2}\biggr)^m
\biggl(\frac{b^2}{k^2}\biggr)^n,
\end{split}
\label{Pcalc}
\end{equation}
where 
\[
\phi_n^{-1} = \frac{2^{\alpha+\beta+1}}{2n+\alpha+\beta+1}
\frac{\Gamma (n+\alpha+1)\Gamma (n+\beta +1)}
{\Gamma (n+1)\Gamma (n+\alpha+\beta+1)},
\]
$x=\cos 2\varphi$, $y=\cos 2\theta$, $a=\sin\varphi \sin \theta$,
$b=\cos \varphi \cos \theta$, $k=(r^{1/2}+r^{-1/2})/2$. The series
(\ref{Pcalc}) is positive for $-1 \leq x$, $y \leq 1$, $0 \leq r <1$ and
$\alpha$, $\beta > -1$.

A Markov process\index{Markov, A. A.!Markov process} analogy to the
Jacobi--Poisson kernel\index{Jacobi, C. G. J.!JacobiPoisson kernel@Jacobi--Poisson kernel|)} is
when the eigenvalues
$c_n(t) = \exp\{-nt\}$.
Following \cite{B1954}\index{Bochner, S.|)} let $\widetilde{X}(t)=X\bigl (Z(t)\bigr )$,
where $\{Z(t), t \geq 0\}$ is a L{\'e}vy
process\index{Levy, P.@L\'evy, P.!L\'evy process} with Laplace
transform\index{Laplace, P.-S.!Laplace transform}
\begin{eqnarray}
&&\mathbb{E}\left [e^{-\lambda Z(t)}\right ]
\nonumber \\
&&~~~~=
\exp\Bigl \{-t \Bigl [\sqrt{2\lambda + (\theta - 1)^2/4} - 
\sqrt{(\theta -1)^2/4}\Bigr ]\Bigr \}\label{SLT}\\
&&~~~~=
\exp \Bigl \{-\frac{t}{\sqrt{2\pi}}\int_0^\infty 
\frac{e^{-x(\theta -1)^2/8}}{x^{3/2}}\bigl (1-e^{-x\lambda}\bigr ) \:dx\Bigr \}.
\nonumber
\end{eqnarray}
$\{Z(t), t \geq 0\}$ is a tilted\index{tilt} positive \Index{stable process}
with index $\frac{1}{2}$ such that $Z(t)$ has an Inverse Gaussian
density\index{Gauss, J. C. F.!inverse Gaussian distribution}
\[
IG\Bigl (\frac{2t}{|\theta - 1|},t^2 \Bigr),\ \theta \ne 1;
\]
that is,
\begin{equation}
\frac{t}{\sqrt{2\pi z^3}}\exp \Bigl \{-\frac{1}{2z}
\Bigl (\frac{|\theta - 1|}{2}z-t\Bigr )^2 \Bigr \},\:z > 0.
\label{igdist}
\end{equation}
The usual \Index{stable density} is obtained when $\theta = 1$ and
(\ref{igdist}) is a tilted density in the sense that it is proportional to
$\exp \bigl \{-z(\theta - 1)^2/8\bigr \}$ times the stable density.
See \cite{F1971} XIII,\index{Feller, W.} \S11, Problem 5 for an early
derivation. $Z(t)$ is distributed as the \Index{first passage time}
\[
T_t = \inf \Bigl \{u>0;B(u) + \frac{|\theta - 1|}{2}u = t\Bigr \},
\]
where $\bigl \{B(u), u \geq 0\bigr \}$ is standard Brownian
motion\index{Brown, R.!Brownian motion}.
The eigenvalues of $\widetilde{X}(t)$ are
\begin{eqnarray}
\widetilde{c}_n(t) &= &
\mathbb{E}\Bigl [\exp \Bigl \{-\frac{1}{2}n(n+\theta -1)Z(t)\Bigr \}\Bigr ]
\nonumber \\
&=&
\exp \Bigl \{-t \Bigl [\sqrt{n(n+\theta-1)+(\theta-1)^2/4} - \sqrt{(\theta -1)^2/4}\Bigr ]\Bigr \}
\nonumber \\
&=&\exp \Bigl \{-t \Bigl [n+(\theta-1)/2 - |\theta -1|/2\Bigr ]\Bigr \}
\nonumber \\
&=&
\begin{cases}
\exp\{- nt\} \text{~~if $\theta \geq 1$},\\
\exp\{-nt\}\times \exp\{t(1-\theta)\} \text{~~if $\theta < 1$.}
\end{cases}
\label{evcalc}
\end{eqnarray}
The process $\{\widetilde{X}(t), t \geq 0\}$ is a
\index{diffusion!jump diffusion}jump diffusion process, discontinuous at the
jumps of $\{Z(t), t \geq 0\}$. Jump sizes increase as $\theta$ decreases.
If $\theta < 1$ then for $n \geq 1$
\[
\mathbb{E}\Bigl [\exp \Bigl \{-\frac{1}{2}n(n+\theta -1)Z(t)\Bigr \}\Bigr ]
= \exp\{- nt\}\times \exp\{t(1-\theta)\},
\]
so subordination\index{subordination|(} does not directly produce eigenvalues
$e^{-nt}$.  Let $\widetilde{f}(x,\allowbreak y;\allowbreak t)$ be the
transition density\index{transition density|(} of $\widetilde{X}(t)$, then the
transition density with eigenvalues $\exp\{-nt\}$, $n \geq 0$ is
\[
e^{-t(1 - \theta)}\widetilde{f}(x,y;t)
+\bigl (1-e^{-t(1-\theta)}\bigr )f_{\alpha\beta}(y).
\]
The subordinated process with this transition density is $X(\widehat{Z}(t))$,
where $\widehat{Z}(t)$ is a similar process to $Z(t)$ but has an extra state
infinity. $Z(t)$ is killed\index{killing} by a jump to infinity at a rate
$(1-\theta)$. Another possible construction does not kill the process
$\widetilde{X}$, but restarts it in a stationary state drawn from the Beta
distribution\index{beta distribution}.  It is convenient to use the notation
that a process $\{Z^\circ (t), t \geq 0\}$ is $\{Z(t), t \geq 0\}$ if
$\theta \geq 1$, or $\{\widehat{Z}(t), t \geq 0\}$ if $0 < \theta < 1$, and use
the single notation $\{X(Z^\circ(t)), t \geq 0\}$ for the subordinated
process.  The transition density (\ref{general:4}), where $c_n(t)$ has the
general form
\begin{equation*}
\exp\left \{-t\int_0^\infty \frac{\bigl (1 - e^{-ny}\bigr )}{y}H(dy)\right \},
\end{equation*}
can then be obtained by a composition of subordinators from the Jacobi
diffusion with any $\alpha$,
$\beta >0$.

There is a question as to which processes with transition densities
(\ref{general:4}) and eigenvalues $c_n(t)$ described by (\ref{ezero}) are
subordinated Jacobi diffusion
processes\index{Jacobi, C. G. J.!Jacobi diffusion|(}. We briefly consider
this question.  Substituting
\[
R_n^{(\alpha ,\beta)}(y) = \>
_2F_1(-n,n+\theta - 1;\beta;1-y)
\]
in the eigenvalue expression (\ref{ezero}),
\begin{eqnarray*}
&&\int_0^{1-}\>\frac{1 - R^{(\alpha ,\beta)}_n(y)}{1-y}\>\nu(dy)\\
&&~~= -c\sum_{k=1}^n
\frac{(-n)_{(k)}(n+\theta - 1)_{(k)}}
{\beta_{(k)}}
\frac{\mu_{k-1}}{k!}\\
&&~~= -c\sum_{k=1}^n
\frac{\prod_{j=0}^{k-1}\bigl (-n(n+\theta-1) + j(j+\theta - 1)\bigr )}
{\beta_{(k)}}
\frac{\mu_{k-1}}{k!},
\end{eqnarray*}
where  $\int_0^{1-}(1-y)^k\>\nu(dy) = c\mu_k$.
The \Index{generator} corresponding to a process with these eigenvalues is
\[
\widehat{{\cal L}}  
= c\sum_{k=1}^\infty
\frac{\prod_{j=0}^{k-1}\bigl (2{\cal L} + j(j+\theta - 1)\bigr )}
{\beta_{(k)}}
\frac{\mu_{k-1}}{k!},
\]
where ${\cal L}$ is the Jacobi diffusion process generator (\ref{generator:0}).
The structure of the class of stochastic processes with the generator 
$\widehat{{\cal L}}$ needs to be understood better. It includes all subordinated Jacobi diffusion processes, but it seems to be a bigger class. A process with generator $\widehat{{\cal L}}$ is a subordinated Jacobi diffusion process if and only if the first derivative of
\begin{equation}
-\sum_{k=1}^\infty
\frac{\prod_{j=0}^{k-1}\bigl (-2\lambda + j(j+\theta - 1)\bigr )}
{\beta_{(k)}}
\frac{\mu_{k-1}}{k!}
\label{CMcond}
\end{equation}
is a completely monotone\index{complete monotonicity} function of $\lambda$. Factorizing 
\[
-2\lambda + j(j+\theta -1) = (j+r_1(\lambda))(j+r_2(\lambda)),
\]
where
$r_1(\lambda),r_2(\lambda)$ are
\[
(\theta -1 )/2 \pm \sqrt{2\lambda + (\theta-1)^2/4},
\]
(\ref{CMcond}) is equal to
\begin{equation}
- \int_0^{1-}\Bigl [ \:
_2F_1(r_1(\lambda),r_2(\lambda);\beta;1-y) -1
\:\Bigr ](1-y)^{-1}\nu (dy).
\end{equation}
\section{Subordinated coalescent process}
Subordinating the Jacobi diffusion
process\index{Jacobi, C. G. J.!Jacobi diffusion|)} $\{X(t),t \geq 0\}$ leads
to subordinating the coalescent\index{coalescence} dual\index{duality}
process, which we investigate in this section.  A subordinated process
$\{\widetilde{X}(t) = X(Z(t)), t \geq 0\}$ has a similar form for the
transition density\index{transition density|)} as (\ref{dual:0}), with
$q^\theta_l (t)$ replaced by $\mathbb{E}(q^\theta_l(Z(t))$, which are
transition functions\index{transition function} of the
subordinated \Index{death process} $A^\theta(Z(t))$. The subordinated process
comes from subordinating\index{subordination|)} the forest of non-mutant
lineages\index{lineage} in a \Index{coalescent tree}.

If $\widetilde{A}^\theta (t) = A^\theta (Z^\circ(t))$, with $Z^\circ (t)$ defined in the last section, we will show that
the probability distribution of $\widetilde{A}^\theta (t)$, $\theta > 0$ is
\begin{equation}
\binom{2k + \theta -1}k\Bigl (\frac{z}{1+z}\Bigr )^k\Bigl (\frac{1}{1+z}\Bigr )^{k+\theta}(1-z),
\label{subdistn}
\end{equation}
for $k\in \mathbb{Z}_+$, where $z=e^{-t}$.  The distribution (\ref{subdistn})
is the distribution of the number of edges in a
time-subordinated \Index{forest}.  Note that if $0 < \theta < 1$ we still
invoke a \Index{subordinator} with a possible jump to infinity at rate
$1-\theta$, so
\[
\mathbb{E}\bigl [q_k^{\theta}(Z^\circ(t)) \bigr ]
= e^{-(1-\theta)t}
\mathbb{E}\bigl [q_k^{\theta}(Z(t)) \bigr ]
+ 
(1 - e^{-(1-\theta)t})\delta_{k0},
\]
because $q_k^\theta(\infty) = \delta_{k0}$.  Although $\theta$ is greater than
zero in (\ref{subdistn}), it is interesting to consider the subordinated
Kingman coalescent\index{Kingman, J. F. C.!Kingman coalescent} with no
mutation\index{mutation|(}. Then $A^0(t) \geq 1$, and
\[
\mathbb{E}\bigl [q_k^{0}(Z^\circ(t)) \bigr ]
= 
e^{-t}\mathbb{E}\bigl [q_k^{0}(Z(t)) \bigr ]
+ 
(1 - e^{-t})\delta_{k1},
\]
because a jump to infinity is made at rate 1, and $q_k^\circ(\infty) = \delta_{k1}$.
The distribution of $\widetilde{A}^0 (t)$ is then, for $k \geq 1$,
\begin{equation}
\binom{2k -1}k\Bigl (\frac{z}{1+z}\Bigr )^k\Bigl (\frac{1}{1+z}\Bigr )^{k}(1-z) + \delta_{k1}(1-z).
\label{csubdistn}
\end{equation}
%% Identity
%%\[
%%\binom{2k-1}k = 2^{2k-1}(-1)^k\left (-\frac{1}{2}\right )_{[k]}/k!
%% \]
The proof of (\ref{subdistn}) ($\theta > 0$) and (\ref{csubdistn}) (with $\theta = 0$) follows directly from the expansion (\ref{trQ}).
\begin{eqnarray}
\mathbb{E}\bigl [q_k^{\theta}(Z^\circ(t)) \bigr ]
&=&\sum_{j=k}^\infty z^j
(-1)^{j-k}
\frac{ (2j+ \theta - 1)(k+\theta)_{(j-1)}}{k!(j-k)!}
\nonumber \\
&=&
\frac{\Gamma (2k+\theta)}{k!\Gamma (k + \theta)}z^k
\nonumber \\
&&~\times\:
\biggl \{1 + \sum_{j=1}^\infty (-1)^j(2j + 2k + \theta -1)
\frac{(2k+\theta)_{(j-1)}}{j!}z^j\biggr \}
\nonumber \\
&=&
\frac{\Gamma (2k+\theta)}{k!\Gamma (k + \theta)}z^k(1-z)(1+z)^{-(2k+\theta)}
\nonumber \\
&=&
\binom{2k + \theta -1}k\Bigl (\frac{z}{1+z}\Bigr )^k\Bigl (\frac{1}{1+z}\Bigr )^{k+\theta}(1-z).
\label{subproof}
\end{eqnarray}
Effectively, in the expansion (\ref{trQ}) of $q_k^\theta(t)$, terms
$\rho_j(t) = \exp\{-\frac{1}{2}j(j+\theta -1)t\}$ are replaced by 
$z^j = \exp\{-jt\}$.
The third line of (\ref{subproof}) follows from the identity, with $|z| < 1$ 
and $\alpha = 2k + \theta$, that
\[
(1-z)(1+z)^{-\alpha} = 1 + \sum_{j=1}^\infty (-1)^j(2j+\alpha -1)\frac{\alpha_{(j-1)}}{j!}z^j,
\]
proved by equating coefficients of $z^j$ on both sides. Of course, for any $|z| < 1$, since (\ref{subdistn}) is a probability distribution,
\begin{equation}
\sum_{k=0}^\infty
\binom{2k + \theta -1}k\Bigl (\frac{z}{1+z}\Bigr )^k\Bigl (\frac{1}{1+z}\Bigr )^{k+\theta}(1-z) = 1.
\label{subdistn:2}
\end{equation}
The probability generating function\index{probability generating function (pgf)|(}
of (\ref{subdistn}) is
\begin{equation}
G_{\widetilde{A}^\theta (t)}(s) = \Bigl (\frac{1-4pqs}{1-4pq}\Bigr )^{-\frac{1}{2}}\Bigl (\frac{1 - \sqrt{1-4pqs}}{2ps}\Bigr )^{\theta - 1},\  \theta > 0,
\label{pgfsub}
\end{equation}
where $p = e^{-t}/(1+e^{-t})$ and $q=1/(1 + e^{-t})$. The calculation needed to show (\ref{pgfsub}) comes from the identity
\begin{equation}
\sum_{k=0}^\infty
\binom{2k + \theta -1}k w^k = 2^{\theta - 1}
\frac{\bigl (1+\sqrt{1 - 4w}\bigr )^{-(\theta - 1)}}{\sqrt {1-4w}},
\label{subdistn:3}
\end{equation}
which is found by substituting 
\[
w = \frac{z}{(1+z)^2}\text{~or~}z = \frac{1 - \sqrt{1-4w}}{1 + \sqrt{1-4w}}
\]
in (\ref{subdistn:2}), then setting 
\[
w = \frac{sz}{(1+z)^2}
\]
in (\ref{subdistn:3}).  The calculations used in obtaining the distribution
and probability generating function are the same as those used in obtaining
the formula (\ref{complexrep}) in Griffiths \cite{G2006}.  There is a
connection with a simple random walk\index{random walk (RW)|(} on $\mathbb{Z}$ with
transitions $j \to j+1$ with probability $p$ and $j \to j-1$ with probability
$q = 1 - p$, when $q \geq p$.  Let the number of steps to hit $-\theta$,
starting from $0$, be $\xi$.  Then $\xi$ has a probability generating function
of
\[
H(s) = \Bigl (\>\frac{1 - \sqrt{1-4pqs^2}}{2ps} \>\Bigr )^{\theta},
\] 
and $\frac{1}{2}(\xi+\theta)$ has a probability generating function
\[
K(s) = \Bigl (\>\frac{1 - \sqrt{1-4pqs}}{2p} \>\Bigr )^{\theta}.
\]
$\widetilde{A}^\theta(t) + \theta$ has the same distribution as the
size-biased\index{size biasing} distribution of $\frac{1}{2}(\xi+\theta)$,
with probability generating function
\[
G_{\widetilde{A}^\theta(t)}(s) = \frac{sK^\prime(s)}{K^\prime (1)},
\]
identical to (\ref{pgfsub}).  In the random walk\index{random walk (RW)|)}
interpretation $\theta$ is assumed to be an integer; however $H(s)$ is
infinitely divisible\index{infinite divisibility}, so we use the same
description for all $\theta > 0$.  Another interpretation is that $K(s)$ is
the probability generating function\index{probability generating function (pgf)|)}
of the total number of progeny in a Galton--Watson branching
process\index{Galton, F.!Galton--Watson process} with
geometric offspring distribution $qp^k$, $k \in \mathbb{Z}_+$, and \Index{extinction}
probability 1, beginning with $\theta$
individuals. See \cite{F1968}\index{Feller, W.} Sections
X.13 and XII.5 for details of the random walk and branching process descriptions.
An analogous calculation to (\ref{subproof}) which is included in Theorem 2.1
of \cite{G2006} is that
\begin{eqnarray}
&&\mathbb{P}\Bigl (\widetilde{A}^\theta (s+t) = j\Bigm| \widetilde{A}^\theta (s) = i\Bigr )
\nonumber \\
&&~~~~= \binom ij\frac{\Gamma (i+\theta) \Gamma (2j + \theta)}
{\Gamma (j+\theta)\Gamma (i+j+\theta )}
z^j(1-z)
\nonumber \\
&&~~~~~~~~~~~~~~~~~\times\:_2F_1(-i + j + 1, 2j + \theta; i + j + \theta;z),
\label{subproof:1}
\end{eqnarray}
where $z = e^{-t}$.
The jump rate from $i \to j$ found from (\ref{subproof:1}) is
\begin{eqnarray}
&&\binom ij\frac{\Gamma (i+\theta) \Gamma (2j + \theta)}
{\Gamma (j+\theta)\Gamma (i+j+\theta )}
\:_2F_1(-i + j + 1, 2j + \theta; i + j + \theta;1),
\nonumber \\
&&~~~=\binom ij B(j+\theta, i-j)^{-1}\int_0^1 x^{2j + \theta - 1}(1-x)^{2(i-j) -2}\>dx\nonumber \\
&&~~~~~=
\begin{cases}
\binom ij\frac{
\Gamma (2i-2j -1)\Gamma (2j + \theta) \Gamma (i+ \theta)
}
{
\Gamma (i-j) \Gamma (j + \theta) \Gamma (2i + \theta -1)
} & \text{if $j= i - 1$, $i - 2$, \ldots,}\\
\frac{\Gamma (2j + \theta)}{\Gamma (j + \theta) j!}\Bigl (\frac{1}{2}\Bigr )^{2j + \theta} \text{if $i = \infty$}.
\end{cases}
\label{jumprates}
\end{eqnarray}
Bertoin \cite{B2009a}, \cite{B2009b}\index{Bertoin, J.} studies the
genealogical structure of trees\index{tree} in an infinitely-many-alleles
branching process model. In a limit from a large initial population size with
rare mutations\index{mutation|)} the \Index{genealogy} is described by a
continuous-state branching
process\index{branching process!continuous-state branching process}
in discrete time with an Inverse\break
Gaussian\index{Gauss, J. C. F.!inverse Gaussian distribution} reproduction
law. We expect that there is a fascinating connection with the process
$\{\widetilde{A}^\theta (t), t \geq 0\}$.  A potential class of transition
functions\index{transition function} of Markov
processes\index{Markov, A. A.!Markov process} $\{\widehat{q}^\theta_k(t),
t \geq 0\}$ which are more general than subordinated processes and related to
Bochner's\index{Bochner, S.}
characterization comes from replacing by $\rho_n^\theta(t)$ by $c_n(t)$
described by (\ref{ezero}); however it is not clear that all such potential
transition functions are positive, apart from those derived by
\Index{subordination}.

\end{document}